\documentclass{amsart}

\usepackage{xypic}
\usepackage{amssymb}
\usepackage{dsfont}
\usepackage{hyperref}
\usepackage{graphicx}
\usepackage{subcaption}
\usepackage{enumitem}
\usepackage{pythonhighlight}

\usepackage{comment}
\usepackage{amsmath}
\usepackage{mathtools}
\usepackage{enumitem} 

\usepackage{tikz-cd}

\xyoption{all}

\newif\iffurther
\furthertrue

\numberwithin{equation}{section}
\numberwithin{figure}{section}
\theoremstyle{plain}
\newtheorem{thm}{Theorem}[section] 
\newtheorem*{thm*}{Theorem}

\newtheorem{prop}[thm]{Proposition}
\newtheorem{cor}[thm]{Corollary}
\newtheorem*{cor*}{Corollary}
\theoremstyle{definition}

\newtheorem{defn}[thm]{Definition}

\newtheorem{lem}[thm]{Lemma}

\theoremstyle{remark}
\newtheorem{rem}[thm]{Remark}

\newtheorem{exmpl}[thm]{Example}

\newtheorem*{acknowledgement*}
{Acknowledgement}

\newcommand\suchthat{\;\ifnum\currentgrouptype=16 \middle\fi|\;}

\def\F{\mathcal{F}}

\def\indeg{{\operatorname{indeg}}}
\def\outdeg{{\operatorname{outdeg}}}

\def \F{\mathcal{F}}

\newcommand{\eps}{\varepsilon}

\title{Non-Amenability of Small Overlap Semigroups}
\author{Margaret Shideler}
\address{Department of Mathematics and Statistics, CUNY Hunter College, 695 Park Avenue, New York NY 10065, USA} 
\email{margaret.shideler@hunter.cuny.edu}

\begin{document}

\begin{abstract}
    We prove that every semigroup admitting a presentation satisfying the small overlap condition $C(3)$ is neither left nor right reversible, unless it is the free semigroup on one generator. Consequently, all such semigroups outside this exception are both left and right non-amenable. We also show the optimality of this result, giving examples of $C(2)$ semigroups that are left and right amenable. Our proof makes use of semigroup diagrams, extending the results of Remmers who first introduced the methods as a means of studying small overlap semigroups. We also compare this to the known result that the $C(7)$ condition in groups implies hyperbolicity which then implies non-amenability (excluding the elementary cases). We establish that, for $m\geq 4$, a strongly $C(m)$ semigroup presentation induces an enveloping group presentation satisfying the small cancellation condition $C(2m-6)$, and show that this bound is sharp.     
\end{abstract}

\maketitle

\section{Introduction}

Small overlap conditions are combinatorial conditions on semigroup presentations analogous to the small cancellation conditions common in combinatorial and geometric group theory \cite{Lyndon}. These conditions were first studied by Remmers \cite{Remmers71, Remmers80} 
in the context of the word problem,  establishing $C(3)$ as a sufficient condition for a solvable word problem. More recent works have studied explicit algorithms for computing normal forms and solving the word problem as well as the relationship between small overlap conditions and other properties such as automaticity \cite{K(SOM1), K(SOM2), K(strongly), K(uniform), Mitchell23, Weiss09}. Small overlap conditions are also \textit{generic} among finitely presented monoids and semigroups, in the sense that for fixed $m, k\geq 1$, and $a\geq 2$, the proportion of $a$-generator, $k$-relation monoid presentations of total relation length $n$ satisfying $C(m)$ tends to 1 as $n$ tends to infinity \cite{K(generic)}. 

The driving insight from the study of small cancellation conditions in groups is that imposing strong enough combinatorial restrictions on presentations has implications for the algebraic, geometric, and algorithmic structure of the group. In particular, one common route for establishing properties of small cancellation groups is to use the well-known result that conditions such as $C(7)$ and $C'(\frac{1}{6})$ give rise to hyperbolic groups \cite{Lyndon}. 

In the semigroup context, the geometric and language-theoretic notions of hyperbolicity need not coincide as they do in groups. Thus, although Kambites showed that $C(4)$ semigroups are word-hyperbolic (in the sense of Duncan and Gilman \cite{Duncan}), this need not imply that all the properties which hold in hyperbolic groups hold in $C(4)$ semigroups. For example, although hyperbolic groups are automatic \cite{Epstein}, Kambites showed that $C(4)$ semigroups were asynchronously automatic (a strictly weaker condition) \cite{K(SOM2)}. Weiss then established the following set of conditions as sufficient to produce automatic semigroups: embeddability into a group, strongly $C(3)$, and for every relation $u=v$, if $u$ is a product of $n_1$ pieces and $v$ is a product of $n_2$ pieces, then $n_1+n_2< 7$ \cite{Weiss09}.

This paper considers which small overlap conditions are sufficient to produce another property present in hyperbolic groups (excluding the elementary cases): non-amenability. Amenability is the existence of a translation-invariant finitely additive probability measure defined on all subsets. For finitely generated groups, non-amenability is witnessed by uniform expansion, which, in non-elementary hyperbolic groups, follows from the geometry of the action on the hyperbolic boundary \cite{Gromov}. For semigroups, the relationship between amenability and expansion is more subtle: translations need not be injective, and the weak and strong F\o lner conditions need not coincide as they do in groups \cite{K(amenable)}. We instead consider another necessary condition for amenability, namely reversibility. Specifically, for a semigroup to be left amenable it must be left reversible, meaning any two nonempty right ideals must intersect.

Excluding the case where the semigroup is isomorphic to the free semigroup on one generator, we establish the following: 
\begin{thm*}[\textbf{\ref{thm: C(3) not reversible}}]
     $C(3)$ semigroups are neither left nor right reversible. 
\end{thm*}
Thus, we obtain the following as an immediate consequence:
\begin{cor*}[\textbf{\ref{thm: C(3) non-amenable}}]  $C(3)$ semigroups are left and right non-amenable. 
\end{cor*}
We also demonstrate the optimality of $C(3)$ as a sufficient condition for non-amenability. Example \ref{exmpl C(2) amenable} gives a $C(2)$ semigroup which is left and right amenable (and thus left and right reversible). We also provide an example of a left amenable $C(2)$ semigroup that highlights the failure of the specific construction of ideals used in the proof of Theorem \ref{thm: C(3) not reversible}. 

We obtain this result using semigroup diagrams, first introduced in Remmers' initial work on small overlap conditions and which are analogous to van Kampen diagrams for groups \cite{Remmers71, Remmers80} (see also \cite{Higgins92}).  An introduction to semigroup diagrams, together with the results of Remmers used here, is presented in Section 2. We establish a new result about when edges in a $C(3)$ diagram can be labeled by pieces (Lemma \ref{lem: initial e not piece}). This is then used to prove the existence of subwords which generate principal ideals that cannot contain words beginning with a different generator and thus witness a failure of reversibility.

We close with a contextualization of this result. First, we explain how non-reversibility (and thus non-amenability) can be inferred from Kambites' work on $C(4)$ semigroups \cite{K(SOM1)}. Finally, we investigate the relationship between small overlap conditions on semigroup presentations and small cancellation conditions on their induced enveloping group presentations. Extending a result of Kashintsev \cite{Kashintsev}, we establish that the induced enveloping group presentation of a strongly \(C(m)\) semigroup satisfies the \(C(2m-6)\) condition (Corollary \(\ref{C(m) implies C(2m-6)}\)) and show that this bound is sharp (Proposition \(\ref{C(m) implies C(2m-6) sharp}\)). These results illustrate that the small overlap conditions sufficient to establish non-amenability directly in semigroups are substantially weaker than those guaranteeing that the induced enveloping group presentation satisfies classical small cancellation criteria for hyperbolicity. In particular, our direct argument establishes non-amenability under the \(C(3)\) condition, whereas guaranteeing a \(C(7)\) induced group presentation requires the stronger condition of strongly \(C(7)\); strongly \(C(5)\) together with \(D(2)\); or strongly \(C(4)\) together with \(D(3)\).

\subsection*{Use of AI} No mathematical ideas, results, or proofs were generated using AI. GPT-6 Astra was used for literature search and for checking grammar and clarity in the text and proofs. The author reviewed all changes, verified all citations, and takes full responsibility for the mathematical content of the paper. 

\section{Preliminaries}

For a given set $X$, denote by $\F_X$ the free semigroup on $X.$ Thus, $\F_X$ consists of all finite words in the alphabet $X$ under the operation of concatenation. For an indexing set $I$, and a set of relations $R = \{(u_i, u_{-i})\mid i\in I\}\subseteq \F_X\times \F_X$,
denote the congruence closure of $R$ by $R^c$. Denote the quotient semigroup of $\F_X$ by $R^c$ as $S=\langle X \mid R\rangle.$ Note that in this notation if  $R = \{(u_i, u_{-i})\mid i\in I\}$, $S$ is often notated as $\langle X\mid \{u_i = u_{-i}\mid  i\in I\}\rangle.$ The generators and relators of $S$ will refer to the elements of $\F_X$ and the elements $u_i$ for $\pm i\in I$, respectively.

We shall henceforth assume that $R= R^c$. Note that nothing is lost in this assumption as the set of relators in $R$ is the same as the set of relators in $R^c$ and $\langle X\mid R\rangle = \langle X\mid R^c\rangle.$

For $w, w'\in \F_X$, we denote equality as words in the free semigroup by $w = w'$ and equality as elements of $S$ by $w \equiv w'.$ Then, $w \equiv w'$ if and only if there exists a finite sequence $W_0, W_1, \dots W_k\in \F_X$ such that $W_0= w$, $W_k = w'$, and for $1\leq j\leq k$, there exists $(u_{i_j}, u_{-i_j})\in R$, and $\beta,\gamma\in \F_X$ such that $W_{j-1} = \beta u_{i_j} \gamma$ and $W_j = \beta u_{-i_j}\gamma$. We call such a sequence a derivation of $w\to w'.$ A derivation is minimal if $k$ is minimized.

For $w\in \F_X$, we denote the length of $w$ in terms of generators from $\F_X$ by $|w|$. For the sequence $W_0, W_1,\dots, W_k\in \F_X$ given above, we say that the relation $(u_{i_j}, u_{-i_j})$ was applied to $W_j$ at position $p$ if $|\beta| = p.$

\subsection{Small Overlap Conditions}

\begin{defn}
    For a semigroup $\langle X\mid R\rangle$, a word $p\in \F_X$ is called a \textbf{piece} if there exists relators $u, u'$ such that $u = \alpha p\beta$ and $u' = \alpha'p\beta'$ and $\alpha \neq \alpha'$ or $\beta \neq \beta'.$
\end{defn}

Note that the case where $u = u'$ is not excluded. Thus, $p$ is a piece exactly when $p$ appears as a subword in distinct relators or in distinct positions within the same relator.

\begin{defn}
    Let $m\in \mathbb{N}$. A semigroup presentation $\langle X\mid R\rangle$ is $\mathbf{C(m)}$ if no relator is a product of fewer than $m$ pieces.
\end{defn}

Every semigroup presentation is $C(1)$, and a semigroup presentation is $C(2)$ if and only if no relator is a proper subword of a distinct relator. If $\langle X\mid R\rangle$ is $C(m)$, then $\langle X\mid R\rangle$ is $C'(m)$ for any $1\leq m'\leq m.$ For a thorough account the $C(m)$ condition, see Remmers \cite{Remmers71}. Note that the $C(m)$ conditions are properties of a particular semigroup presentation, rather than of the semigroup itself. We will say that a semigroup is a \textit{$C(m)$ semigroup} if it admits a presentation satisfying $C(m)$.

\subsection{Semigroup Diagrams} 

We now give a description of the machinery behind our main result, semigroup diagrams. These were first introduced by Remmers as a means of proving that $C(3)$ semigroups have a solvable word problem \cite{Remmers71}. For a more thorough introduction see \cite{Remmers71} or \cite{Higgins92}.

\begin{defn} {\cite[Definition~1.7.1]{Higgins92}} Let $X$ be a set. A \textbf{semigroup diagram} $M$ on $X$ is a planar connected digraph $\Gamma$ embedded in $\mathbb{R}^2$, together with a labeling function which assigns to each oriented edge of $\Gamma$ an element of $X$ such that
\begin{enumerate}
    \item Each component of $S^2 - \Gamma$ is \textit{two-sided}, meaning that it is bounded by a counterclockwise cycle of the form $\pi\rho^{-1}$, where $\pi$ and $\rho$ are non-empty directed paths where $\rho^{-1}$ denotes the reverse of $\rho$.
    \item $\Gamma$ has exactly one source, $O$, and one sink, $I$, each of which lies on the boundary of the unbounded component of $S^2-\Gamma.$
\end{enumerate}
\end{defn}

If $\pi \rho^{-1}$ is the counterclockwise boundary of the unbounded component of $S^2-\Gamma$, then $\pi$ is referred to as the \textit{left boundary of $M$}, denoted $\alpha_M$, and $\rho $ is referred to as the \textit{right boundary of $M$}, denoted $\omega_M.$ The bounded components of  $S^2-\Gamma$ are referred to as the \textit{regions} of $M$. If $\pi \rho^{-1}$ is the counterclockwise boundary of a region $D$, then $\pi$ is referred to as the \textit{left side of $D$}, denoted $\alpha_D$, and $\rho $ is referred to as the \textit{right side of $D$}, denoted $\omega_D.$ 

For an oriented edge $e$ in $M$ denote by $\overline{e}$ the label of $e$ such that $\overline{e}\in X$. If $\pi$ is a path consisting of consecutive edges $e_1, \dots, e_n$, then $\overline{\pi} =\overline{e_1}\cdots\overline{e_n}\in \F_X.$ The path labeled by the empty word $1$ will consist of a single vertex and no edges.

\begin{defn}
    Let $S = \langle X\mid R\rangle.$ Let $M$ be a semigroup diagram on $X$. Then $M$ is an \textbf{$S$-diagram} if $(\overline{\alpha_D}, \overline{\omega_D})\in R$ for every region $D$ of $M$. Suppose $w, w'\in \F_X$ such that $w \equiv w'$. Then, $M$ is called an $S$-diagram for the derivation of $w \to w'$ if $\overline{\alpha_M}= w$ and  $\overline{\omega_M} = w'.$ $M$ is called minimal if the derivation $w\to w'$ is minimal. 
\end{defn}

Suppose the derivation $w \to w'$ is given by the sequence $w=W_0, W_1, \cdots, W_k = w'$ such that for every $1\leq j\leq k$, $$W_{j-1}  = \beta_ju_{i_j}\gamma_j \quad \text{and} \quad W_j = \beta_j u_{-i_j}\gamma_j$$
and $(u_{i_j}, u_{-i_j})\in R.$ Then the unique $S$-diagram for this derivation can be constructed as follows. If $k = 0$, then $w = w'$ and $M$ is a simple directed arc with the label $w.$ If $k>0$, assume inductively that $N$ is an $S$-diagram for the derivation of $W_0\to W_{k-1}$ given by the sequence $W_0, W_1, \cdots, W_{k-1}.$ Thus, $\overline{\omega_N} = W_{k-1}$ so $\omega_N = \pi\sigma\rho$ such that $\overline{\pi} = \beta_{k}$, $\overline{\sigma} = u_{i_k}$, and $\overline{\rho} = \gamma_k.$ To construct $M$ from $N$, attach a simple directed arc $\tau$ with label $u_{-i_k}$ to the initial and terminal vertices of $\sigma$ such that the cycle $\sigma \tau^{-1}$ is oriented counterclockwise and is the boundary of a region $D$. Then, $\omega_M = \pi\tau\rho$, so $\overline{\omega_M} = \beta_k u_{-i_k}\gamma_k = W_k.$ A proof of the uniqueness of this diagram is given by {\cite[Theorem~1.7.2]{Higgins92}}.

\begin{figure}
    \centering
    \begin{subfigure}{0.45\textwidth}
        \centering
        \includegraphics[width=3.7cm]{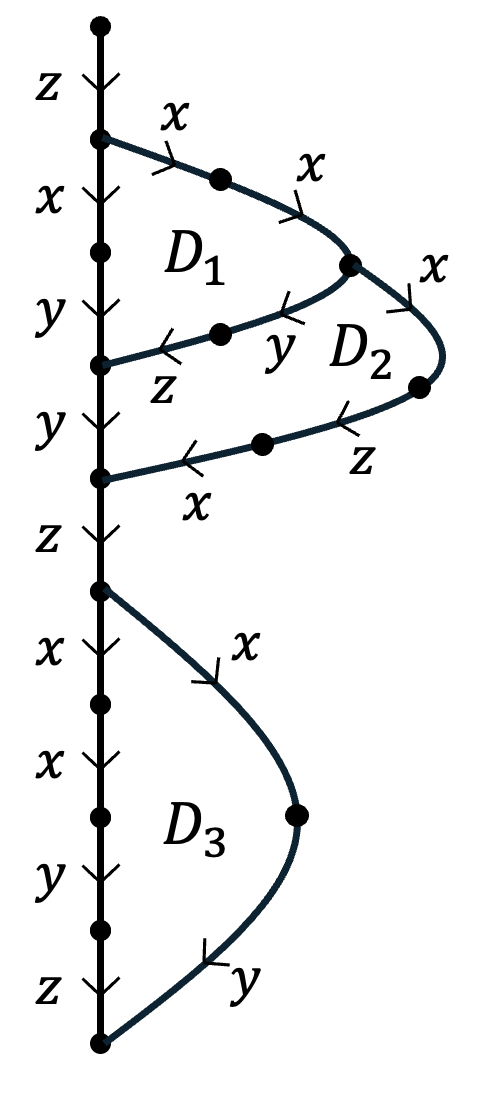}
        \caption{Initial $S$-Diagram after construction}
        \label{fig:exmpl_diagram1}
    \end{subfigure}
    \hfill
    \begin{subfigure}{0.45\textwidth}
        \centering
        \includegraphics[width=4cm]{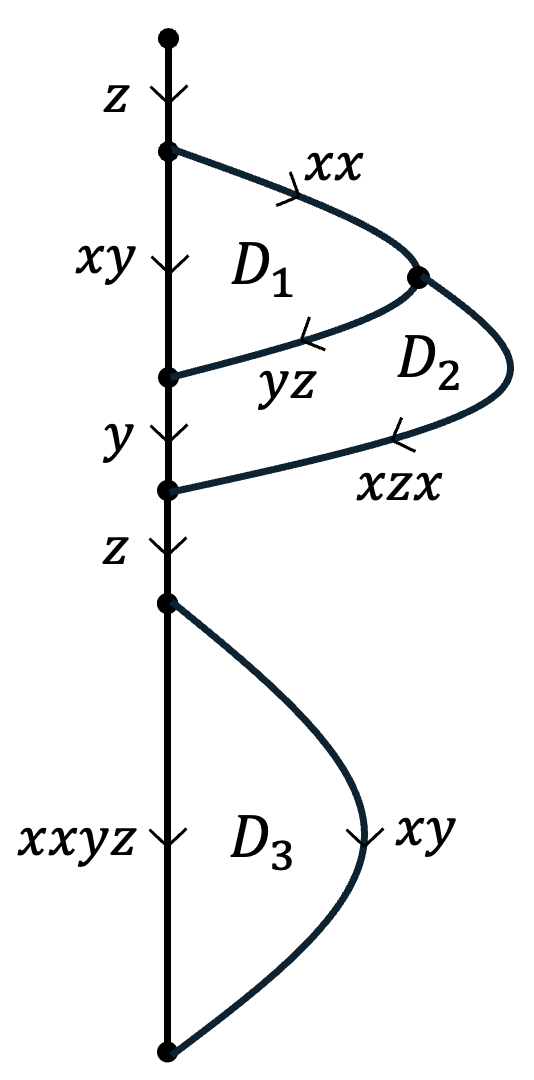}
        \caption{Reduced $S$-Diagram}
        \label{fig:exmpl_diagram2}
    \end{subfigure}
    \caption{$S$-Diagrams for the Derivation $zxyyzxxyz\to zxxxzxzxy$}
\end{figure}

\begin{exmpl}\label{ex:diagram}
    Consider the semigroup $S = \langle x,y,z \mid xy= xxyz, xzx= yzy\rangle.$ We will build the $S$-diagram for the following derivation for $zxyyzxxyz\to zxxxzxzxy:$
    $$z(xy)yzxxyz, z[xx(yz]y)zxxyz, zxx[xzx]z(xxyz), zxxxzxz[xy],$$
    where $(w)$ indicates the relator word to be replaced in the next step and $[w]$ indicates the relator word which has just been replaced.

    Figure \ref{fig:exmpl_diagram1} shows the $S$-diagram for this derivation. Note that the boundaries of the diagram and regions in this diagram are as follows:
    \begin{eqnarray*}
        \text{Diagram}: &\overline{\alpha_S}=zxyyzxxyz \quad &\overline{\omega_S} = zxxxzxzxy \\
       \text{Region } D_1: &\overline{\alpha_{D_1}} = xy  &\overline{\omega_{D_1}} = xxyz\\
    \text{Region } D_2: &\overline{\alpha_{D_2}} = yzy  &\overline{\omega_{D_2}} = xzx\\
       \text{Region } D_3: &\overline{\alpha_{D_3}} = xxyz  &\overline{\omega_{D_3}} = xy
    \end{eqnarray*}
\end{exmpl}

For a vertex $P$ in an $S$-diagram, the \textit{indegree} of $P$, denoted $\indeg (P)$, is the number of edges whose terminal vertex is $P$. The \textit{outdegree} of $P$, denoted $\outdeg(P)$, is the number of edges whose initial vertex is $P$. A vertex $P$ is \textit{superfluous} if $\indeg (P) = \outdeg(P)= 1.$

\begin{defn}
    A semigroup diagram is a \textbf{reduced diagram} if there are no superfluous vertices.
\end{defn}

A reduced diagram $M'$ is obtained from a diagram $M$ in the following manner: Suppose $P$ is a superfluous vertex such that $P$ is the terminal vertex of $e_1$ and the initial vertex of $e_2.$ Denote the initial vertex of $e_1$ and the terminal vertex of $e_2$ by $Q_1$ and $Q_2$, respectively. Replace $e_1, e_2$ and $P$ with an edge $f$ whose initial vertex is $Q_1$ and terminal vertex is $Q_2$. Let $\overline{f} = \overline{e_1}\cdot \overline{e_2}.$ Repeat this process until no superfluous vertices remain. Note that in reduced diagrams, edges are labeled by elements of $\F_X$ rather than by elements of $X$. 

Under this construction, every region of $M$ naturally corresponds to a region of $M'$; corresponding regions have the same left and right side labels. Furthermore, $\overline{\alpha_{M'}}=\overline{\alpha_M}$ and $\overline{\omega_{M'}}=\overline{\omega_M}$. Unless otherwise specified, all diagrams will be assumed to be reduced. 

\begin{exmpl}
    The reduced diagram corresponding to the derivation presented in Example \ref{ex:diagram} is shown in Figure \ref{fig:exmpl_diagram2}.
\end{exmpl}

\subsection{Preliminary Results on Semigroup Diagrams} Suppose that $D$ is a region of a diagram $M$. A \textit{side} of $D$ refers to either $\alpha_D$ or $\omega_D$. The \textit{initial edge} (resp. \textit{terminal edge}) of a side is the edge that shares an initial (resp. terminal) vertex with the side. An \textit{intermediate edge} of a side is any edge which is not an initial or terminal edge. An edge is \textit{interior} if it is not a subpath of $\alpha_M$ or $\omega_M$. A side is \textit{interior} if all its component edges are interior. Two regions $D$ and $E$ share a \textit{common edge} $e$ if $e$ is a subpath of both $\alpha_D$ and $\omega_E$ or $\omega_D$ and $\alpha_E$. 

We recall the following results proved by Remmers in \cite{Remmers71}:

\begin{lem}\label{lem:remmers}
    Let $S = \langle X\mid R\rangle$. Suppose $M$ is an $S$-diagram.

    \begin{enumerate}[label=(\roman*), ref=\thelem(\roman*)]
    \item\label{lem:remmers:common edge piece} {\cite[Lemma~4.2]{Remmers71}} Suppose that $D$ and $E$ are regions with common edge $e$ such that $e$ is the initial edge of $\alpha_D$ but not of $\omega_E$ or the terminal edge of $\alpha_D$ but not of $\omega_E$. Then, $\overline{e}$ is a piece.

    \item\label{lem:remmers:no one edge regions} {\cite[Lemma~4.5]{Remmers71}} Suppose that $\langle X\mid R\rangle$ is $C(2)$, and $M$ is a minimal diagram. Then no interior edge of $M$ is an entire side of a region. 

    \item\label{lem:remmers:three edges max} {\cite[Lemma~4.8a]{Remmers71}} Suppose that $\langle X\mid R\rangle$ is $C(3)$. Then no side of any region has more than 3 edges.

    \item\label{lem:remmers:three edges interior} {\cite[Lemma~4.8b]{Remmers71}} Suppose that $\langle X\mid R\rangle$ is $C(3)$. If two consecutive edges on a side of region are interior, then the entire side of the region must be interior.
    \end{enumerate}
\end{lem}

\subsection{(Non)-Amenability of Semigroups} Finally, we recall the definition of amenability for semigroups as well as known necessary and sufficient conditions for semigroup amenability. For a more in-depth discussion, see \cite{K(amenable)}. 

\begin{defn}
    A semigroup $S$ is \textbf{left amenable}  if it admits a finitely additive probability measure $\mu$, defined on all the subsets of $S$, which is left invariant, in the sense that $\mu(\{s\in S\mid as\in X\}) = \mu(X)$ for all $X\subseteq S$ and $a\in S$. 
\end{defn}

\begin{defn}
    A semigroup $S$ satisfies the \textbf{(left) F\o lner condition} (FC) if for every finite subset $H\subseteq S$ and every $\eps >0$, there is a finite non-empty subset $F$ of $S$ with $|sF\setminus F|\leq \eps |F|$ for all $s \in H$. 

    A semigroup $S$ satisfies the \textbf{(left) strong F\o lner condition} (SFC) if for every finite subset $H\subseteq S$ and every $\eps >0$, there is a finite non-empty subset $F$ of $S$ with $|F\setminus sF|\leq \eps |F|$ for all $s \in H$. 
\end{defn}

\begin{defn}
    A semigroup $S$ is \textbf{left reversible} if every pair of right ideals intersect. Thus, for every $a,b\in S$, $$aS\cap bS \neq \varnothing.$$
\end{defn} 

The above conditions are related in the following way (see \cite{K(amenable)}).
\[
\begin{array}{ccccc}
\text{left SFC} & \Longrightarrow & \text{left amenable}
                & \Longrightarrow & \text{left FC} \\
                & & \Downarrow & & \\
                & & \text{left reversible} & &
\end{array}
\]
Notably, none of these implications hold in the converse. 

 In particular, we will take advantage of relationship between amenability and reversibility to show that $C(3)$ semigroups are left and right non-amenable after establishing they neither left nor right reversible.

\section{Main Results}

We now establish one new result on semigroup diagrams for $C(3)$ semigroups which will be used to prove the main theorem. 

\begin{lem}\label{lem: initial e not piece}
    Suppose $S = \langle X \mid R\rangle$ is $C(3)$ and that $M$ is a minimal $S$-diagram. Let $e$ be an edge of the left (resp. right) boundary of $M$ that is also an initial edge of the left (resp. right) side of some region $D$. Then, $\overline{e}$ must not be a piece.
\end{lem}

\begin{figure}
    \centering
    \begin{subfigure}{0.45\textwidth}
        \centering
        \includegraphics[width=3cm]{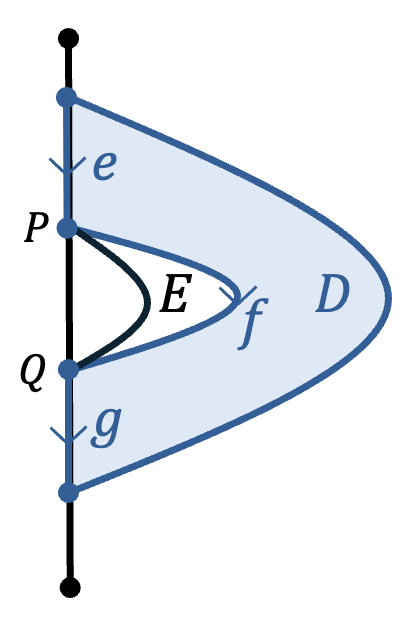}
        \caption{Case where $\alpha_D$ has three edges}
        \label{fig:three_edges}
    \end{subfigure}
    \hfill
    \begin{subfigure}{0.45\textwidth}
        \centering
        \includegraphics[width=2.3cm]{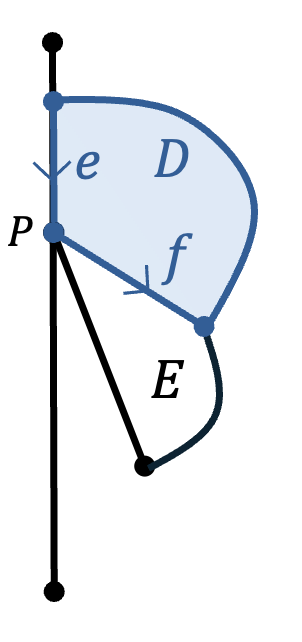}
        \caption{Case where $\alpha_D$ has two edges}
        \label{fig:two_edges}
    \end{subfigure}
    \caption{$S$-Diagrams where the initial edge of $\alpha_D$ is a subpath of $\alpha_M$}
\end{figure}

\begin{proof}
    Without loss of generality, suppose $e$ is a subpath of $\alpha_M$ and is the initial edge of $\alpha_D$. Denote the terminal vertex of $e$ by $P$. Since $\langle X\mid R\rangle$ is $C(3)$, by Lemma \ref{lem:remmers:three edges max}, $\alpha_D$ must have no more than three edges. 
    
    Suppose that $\alpha_D$ has three edges. Denote by $f$ and $g$ the intermediate and terminal edges of $\alpha_D$, respectively. Since $P$ is not superfluous,  $f$ must be interior. And since $\outdeg(P)>1$, $f$ must be the initial edge of the side $\omega_E$ of some region $E$. Denote the terminal vertex of $f$ by $Q$. By Lemma \ref{lem:remmers:three edges interior}, $g$ cannot be interior as then $e$ would be interior. Thus, $\indeg(Q)>1$, so $f$ must also be the terminal edge of $\omega_E$. Thus, $\omega_E$ is a side of a region composed of one interior edge, contradicting \ref{lem:remmers:no one edge regions}, see Figure \ref{fig:three_edges}. Therefore, $\alpha_D$ must have one or two edges.

    Suppose $\alpha_D$ has one edge, so $\alpha_D = e.$ Then, if $\overline{e}$ is a piece, $\overline{\alpha_D}$ is a relator which is a piece, contradicting that $\langle X\mid R\rangle$ is $C(2).$

    Finally, suppose $\alpha_D$ has two edges. Denote the terminal edge of $\alpha_D$ by $f$, see Figure \ref{fig:two_edges}. Since $P$ is not superfluous, $\outdeg(P) >1$, so $f$ is the initial edge of $\omega_E$ for some region $E$. By Lemma \ref{lem:remmers:no one edge regions}, $f$ must not be the terminal edge of $\omega_E.$ Thus, by Lemma \ref{lem:remmers:common edge piece}, $\overline{f}$ is a piece. Then, if $\overline{e}$ is a piece, $\overline{\alpha_D} = \overline{e}\cdot \overline{f}$ is a decomposition of a relator into 2 pieces, contradicting that $\langle X\mid R\rangle$ is $C(3).$
\end{proof}

Symmetrically, we obtain the following: 
\begin{cor}\label{cor: terminal e not piece}
    Suppose $S = \langle X \mid R\rangle$ is $C(3)$ and that $M$ is an $S$-diagram. Let $e$ be an edge of the left (resp. right) boundary of $M$ that is also an terminal edge of the left (resp. right) side of some region $D$. Then, $\overline{e}$ must not be a piece.
\end{cor}

\begin{rem}\label{rem: assumptions}
    Before establishing the main result about $C(3)$ semigroups, we must consider the case where the presentation results in a semigroup isomorphic to a free semigroup. In particular, a free semigroup is left and right amenable if and only if it is of rank one. 

    Thus, we will henceforth make the following assumptions when considering a $C(3)$ semigroup $S = \langle X\mid R\rangle.$
    \begin{enumerate}
        \item $R$ contains at least one relation.
        \item $X$ contains at least two generators. Note that this assumption is natural as semigroup with a relation of the form $x^a = x^b$ for $x\in X$, $a,b\in \mathbb{N}$ would fail to be $C(2).$
        \item At least one relator has a length greater than 1. Note that if all relations of $S$ are of length 1, then $S$ is isomorphic to the free semigroup on a subset of $X$.
    \end{enumerate}
\end{rem}

\begin{thm}\label{thm: C(3) not reversible}
     $C(3)$ semigroups are not left reversible. Symmetrically, $C(3)$ semigroups are not right reversible.
\end{thm}

\begin{proof} Let $S$ be a semigroup with a $C(3)$ presentation $ \langle X\mid R\rangle$.

\smallskip
    Suppose that there exists some relator $u = x\beta$ where $x\in X$ and $\beta$ is not a piece. Then, consider a minimal $S$-diagram $M$ with $\overline{\alpha_M} = \beta W$ for some $W\in F_X.$ Suppose that there is a region $D$ whose initial edge $e$ of $\alpha_D$ is the initial edge of $\alpha_M$ (i.e., $\overline{e}$ is the prefix of $\beta W$). 
    
    Suppose that $|\overline{e}|\leq |\beta|$, so $\overline{e}$ is the prefix of $\beta$. Then, we have decompositions of relators $u = x\cdot\overline{e}\cdot\beta'$ and $\overline{\alpha_D} = \overline{e} \cdot \gamma$ for some $\beta', \gamma\in F_X$. Thus, $\overline{e}$ is a piece contradicting Lemma \ref{lem: initial e not piece}. Suppose instead that $|\overline{e}|>|\beta|$, so $\beta$ is the prefix of $\overline{e}$. Then, we have decompositions of relators $u = x\cdot\beta$ and $\overline{\alpha_D} = \beta \cdot\delta$ for some $\delta\in F_X$. Thus, $\beta$ is a piece, contradicting our assumption.
    
    Thus, the initial edge of $\alpha_M$ must not be an edge in the side of any region. In particular, suppose that $y\in X$ is not the first generator of $\beta$. Recall by the assumptions in Remark \ref{rem: assumptions} such a generator exists. Then, $\omega_M$ cannot be of the form $yW'$ for any $W'\in F_X$ as this would require a region to have the the initial edge of $\alpha_M$ as part of its boundary. Thus, $\beta S\cap yS = \varnothing$, so $S$ is not left reversible.

\smallskip
    Suppose instead that no relator can be expressed in the form $u = x\beta$ where $x\in X$ and $\beta$ is not a piece. As $\langle X\mid R\rangle$ is $C(3)$, if $\beta$ is a piece, then $x$ must not be a piece. Thus, the first generators of all relators must not be pieces. By assumption, there exists at least one relator, $u$, of length greater than 1. Let $x$ be the initial generator of $u$ and $y$ be any other generator which appears in $u$. Then, $y$ does not appear as the first generator of any relator since then it would be a piece. Thus, no relation can be applied at the 0th position of $yW$ for any $W\in F_X$. Thus, $yS\cap xS = \varnothing,$ so $S$ is not left reversible.

\smallskip
    The symmetrical result for non-right reversibility is obtained by considering a decomposition of a relator $u = \beta x$ where $x\in X$ and $\beta$ not a piece and using Corollary \ref{cor: terminal e not piece}.
\end{proof}

Thus, using the fact that reversibility is a necessary condition for amenability, we establish the following:

\begin{cor}\label{thm: C(3) non-amenable}
     $C(3)$ semigroups are left and right non-amenable.
\end{cor}

We now show the above result is optimal and provide an example of a $C(2)$ semigroup which is left and right amenable (and thus left and right reversible). 

\begin{exmpl}\label{exmpl C(2) amenable}
    Consider  $S = \langle a, b \mid  ab = ba\rangle.$ This presentation is $C(2)$ but not $C(3)$ as neither relator is a subword of the other but both relators are products of two pieces. Since $S$ is commutative, $S$ is left and right amenable {\cite[Theorem~4]{Argabright}}. Note that $C(2)$ is strictly the strongest condition which can hold in a commutative semigroup.
\end{exmpl}

We also provide an example of a left amenable $C(2)$ semigroup that witnesses the failure of the argument used in Theorem \ref{thm: C(3) not reversible}. 

\begin{exmpl}\label{exmpl2 C(2) amenable}
    Consider $S = \langle a, b \mid abb = bab\rangle.$ $S$ is $C(2)$ as neither relator is a subword of the other, and $S$ is not $C(3)$ as $abb = (ab)\cdot b$ is a product of 2 pieces. 

    To show $S$ is left amenable, we will show that $S$ satisfies the left strong F\o lner condition. For $N\in \mathbb{N}$, let $F_N = \{b^i(ab)^j \mid 1\leq i,j\leq N\}.$ We can use $aab\to bab$ as a string rewriting system for $S$. This is confluent as $aab$ contains no self-overlaps. Thus, all elements in $F_N$ are irreducible so $|F_N| = N^2.$ Since $S$ is finitely generated, it suffices to show that for every $\eps >0$, there exists $N\in \mathbb{N}$, such that $|F_N\setminus aF_N|, |F_N\setminus bF_N| \leq \eps |F_N|$. 

    Consider $aF_N$. Since $ab^i(ab)^j\equiv  b^{i-1}(ab)^{j+1}$, $$F_N \setminus aF_N = \{b^N(ab)^j\mid 2\leq j\leq N\}\cup \{b^i(ab) \mid 1\leq i\leq N\},$$
    so $|F_N\setminus aF_N| = 2N-1.$

    Consider $bF_N$. $$F_N\setminus bF_N =\{b(ab)^j \mid 1\leq j\leq N\},$$
    so $|F_N\setminus bF_N| = N.$

    Thus, $S$ satisfies the left SFC and is left amenable. Note that this example highlights the failure of the argument used in Theorem \ref{thm: C(3) not reversible} as there is a relator $abb = a\cdot bb$ where $bb$ is not a piece.
\end{exmpl}

We conclude this section by situating this result within other known results. Beyond Remmers' 1971 result about the word problem for $C(3)$ semigroups, other breakthroughs about properties of small overlap semigroups have come from Kambites' work on $C(4)$ semigroups: $C(4)$ semigroups have linear time solvable word problem, and a regular language of linear-time computable normal forms, and are also rational, asynchronous automatic and word hyperbolic, see \cite{K(SOM1), K(SOM2), K(strongly)}. In particular, in \cite{K(SOM1)}, Kambites presents an algorithm for the word problem of a $C(4)$ presentation. There he gives the following proposition:

\begin{prop}{\cite[Proposition~3]{K(SOM1)}}
    Fix a $C(4)$-presentation for a semigroup $S$. Suppose a word $u$ has an overlap prefix $aXY$ and that $u = aXYu'$. Then $u\equiv v$ if and only if
$v = av'$ where $v \equiv XYu'.$
\end{prop}

In his notation, an overlap prefix $aXY$ is a factorization of the form $$aX_1Y_1'X_2Y_2'\cdots X_nY_n$$ where:
\begin{itemize}
    \item $n\geq 1$;
    \item $aX_1Y_1'X_2Y_2'\cdots X_nY_n$ has no factor of the form $X_0Y_0$, where $X_0$ and $Y_0$ are the maximal piece prefix and middle word, respectively, of some relation word, beginning before the end of the
prefix $a$;
    \item for each $1\leq i\leq n$, $R_i = X_iY_iZ_i$ is a relation word with $X_i$ and $Z_i$ the maximal piece prefix and suffix respectively; and
    \item for each $1\leq i<n$, $Y_i'$ is a proper, non-empty prefix of $Y_i$
\end{itemize}

From this result, one could conclude that $C(4)$ semigroups are not right reversible, for if $a,b$ are distinct generators $(aXY)S \cap bS = \varnothing.$ Theorem \ref{thm: C(3) not reversible} strengthens this result.

\section{Small Cancellation Conditions of Enveloping Groups}

 The group-theoretic analogue of small overlap conditions, called small cancellation conditions, are well-understood \cite{Lyndon}.  In particular, several classical small cancellation conditions are sufficient to establish hyperbolicity and, excluding elementary cases, non-amenability. This raises the question of whether the non-amenability result of Corollary~\ref{thm: C(3) non-amenable} could be obtained by passing to the enveloping group and applying these classical hyperbolicity criteria. To investigate this possibility, we compare the small overlap conditions of a semigroup presentation with the small cancellation conditions of its induced enveloping group presentation.

 Given a semigroup $S$, let $G(S)$ denote its enveloping group with homomorphism $\iota: S\to G(S)$, characterized by the universal property that every homomorphism $f\colon S\to H$ into a group $H$ factors uniquely through $\iota$. If $S$ has presentation $\langle X\mid u_i=v_i,\ i\in I\rangle$, then the \textit{induced group presentation} of $G(S)$ is obtained by taking the cyclic reductions of $u_iv_i^{-1}$ as relators, for each $i\in I$.

We recall the following definitions from small cancellation theory:

\begin{defn}
    Let $G = \langle X\mid R\rangle$ be a group whose relators $r_i\in R$ are cyclically reduced. Let $R^*$ denote the symmetrized closure of $R$, consisting of all cyclic permutations of the elements of $R^{\pm1}$. A word $p$ is called a \textbf{piece} if there exist distinct words $r_1,r_2\in R^*$ such that $p$ is a common initial segment of both $r_1$ and $r_2$.
    
    For $m\in \mathbb{N}$, a group presentation $\langle X\mid R\rangle$ is $\mathbf{C(m)}$ if no word in $R^*$ can be expressed as a product of fewer than $m$ pieces. 

    For $q\geq 3$, a group presentation $\langle X\mid R\rangle$ is $\mathbf{T(q)}$ if, for every sequence $r_1,\ldots,r_k\in R^*$ with $3\leq k<q$ such that $r_{i+1}\neq r_i^{-1}$ for all $1\leq i\leq k-1$ (and $r_1\neq r_k^{-1}$), at least one of the concatenations $r_ir_{i+1}$ (or $r_kr_1$) is freely reduced; that is, there is no cancellation between the last letter of $r_i$ and the first letter of $r_{i+1}$ (or between the last letter of $r_k$ and the first letter of $r_1$).
\end{defn}

A group admitting a finite presentation satisfying $C(7)$, or more generally $C(m)$--$T(q)$ with $m,q\geq 3$ and
\[
\frac{1}{m}+\frac{1}{q}<\frac{1}{2},
\]
is hyperbolic \cite{Lyndon}. If the group is neither finite nor virtually cyclic, it is non-elementary hyperbolic and hence non-amenable.

\medskip

The relationship between the small overlap condition of semigroup presentations and the small cancellation condition of the induced presentation of their enveloping groups was first investigated by Kashintsev \cite{Kashintsev}. We recall the following definitions, notation, and results. 

First, we will restrict this question to the case where all the relators in a semigroup presentation are distinct using the following terminology:

\begin{defn}
    Let $m\in \mathbb{N}$. A semigroup presentation $\langle X\mid R\rangle$ is \textbf{strongly $\mathbf{C(m)}$} if it is $C(m)$ and has no repeated relators.
\end{defn}

Note that without this assumption, a semigroup presentation with relations of the form $u =v$ and $u=w$ may have relators $uv^{-1}, uw^{-1}$ in the induced presentation of its enveloping group. Thus, $u$ could be a piece with respect to the group presentation, but not the semigroup presentation. For a more thorough discussion on the difference between $C(m)$ and strongly $C(m)$ conditions, see \cite{K(strongly)}. 

\begin{defn}
    Let $\langle  X \mid R\rangle$ be a semigroup presentation with $R = \{(u_i, v_i)\mid i\in I\}$. The \textbf{left graph} and \textbf{right graph} of the presentation
    are undirected multigraphs with vertex set $X$, defined as follows:
    \begin{itemize}
        \item The left graph has an edge $e_i$ joining the
        first generators of $u_i$ and $v_i$ for each $i\in I$.
        \item The right graph has an edge $e_i'$ joining the
        last generators of $u_i$ and $v_i$ for each $i\in I$.
    \end{itemize}
    
    For $q\in \mathbb{N}$, the presentation $\langle X\mid R\rangle$ satisfies the $\mathbf{D(q)}$ \textbf{condition} if the left and right graphs contain no cycles of length less than $q$. 
\end{defn}

All semigroup presentations are $D(1)$. A semigroup presentation is $D(2)$ if and only if there is no relation $u_i= v_i$ such that $u_i$ and $v_i$ begin or end with the same generator. A semigroup presentation satisfies the $D(3)$ condition if and only if it satisfies $D(2)$ and, for every pair of distinct generators $x_a,x_b\in X$, there is at most one relation composed of relators that begin with $x_a$ and $x_b$, in either order, and  at most one relation composed of relators that end with $x_a$ and $x_b$, in either order.

We now recall the following result of Kashintsev:

\begin{prop} {\cite[Lemma~3.6]{Kashintsev}} \label{prop: kashintsev}
Let $m\geq 2$, and suppose the semigroup presentation $S = \langle X\mid R\rangle$ is strongly $C(m)$. If this presentation of $S$ is $D(2)$, then the induced presentation of the enveloping group $G(S)$ is $C(2m-2)$. If this presentation of $S$ is $D(3)$, then the induced presentation of the enveloping group $G(S)$ is $C(2m)$.
\end{prop}

We extend this result to the case where the presentation need not satisfy $D(2)$. 

\begin{cor} \label{C(m) implies C(2m-6)} 
    Let $m\geq 4$, and suppose the semigroup presentation $S = \langle X\mid R\rangle$ is strongly $C(m)$. Then the induced presentation of the enveloping group $G(S)$ is $C(2m-6).$
\end{cor}

\begin{proof}
    Let $R = \{(u_{i},v_i) : i\in I\}$. For each $i\in I$, let $u_i = \alpha_i s_i\beta_i$ and $v_i = \alpha_i t_i\beta_i$ such that $\alpha_i$ and $\beta_i$ are the longest shared prefix and suffix, respectively, of $u_i$ and $v_i$. If $\alpha_i$ and $\beta_i$ are non-empty, then they are pieces. Since $\langle X\mid R\rangle$ is $C(m)$, no relator can be the product of fewer than $m$ pieces, and thus no word in $\{s_i, t_i\mid i\in I\}$ can be product of fewer than $m-2$ pieces. In particular, since $m\geq 4$, for all $i\in I$, $s_i$ and $t_i$ must be non-empty words that are themselves not pieces. Note this means the words in $\{s_i, t_i\mid i\in I\}$ are pairwise distinct. Since the $s_i, t_i$ are not pieces, if a pair of them were equal then their corresponding relators would be equal, contradicting that $\langle X\mid R\rangle$ is \textit{strongly} $C(m)$. 

    Let $R' = \{(s_i, t_i): i\in I\}$ and consider $S' = \langle X\mid R'\rangle$. Note that if $\langle X\mid R\rangle$ is $D(2)$ then all the $\alpha_i, \beta_i$ are empty words and $R = R'.$ We have that all relators in $\langle X\mid R'\rangle$ are distinct, and if a word $p$ is a piece with respect to $\langle X\mid R'\rangle$ then it must also be a piece with respect to $\langle X\mid R\rangle$, so $\langle X\mid R'\rangle$ is strongly $C(m-2)$. And since $\alpha_i, \beta_i$ are the maximal shared prefixes and suffixes, each pair $s_i, t_i$ must not share the same first or same last generator. Thus, $\langle X\mid R'\rangle$ is $D(2).$ Therefore, by Proposition \ref{prop: kashintsev}, the induced presentation of $G(S')$ is $C(2(m-2)-2) = C(2m-6).$

    Finally, we compare the induced presentations of $G(S)$ and $G(S')$. For every $i\in I$, $u_iv_i^{-1} = \alpha_i s_i\beta_i\beta_i^{-1}t_i^{-1}\alpha_i^{-1}$, so the cyclic reduction of $u_iv_i^{-1}$ is $s_it_i^{-1}.$ Thus, $S$ and $S'$ induce the same presentation for their enveloping groups. 
\end{proof}

\begin{prop}
    The condition given in Corollary \ref{C(m) implies C(2m-6)} is sharp. For $m\geq 4$, there exists a semigroup presentation that is strongly $C(m)$ whose induced enveloping group presentation is not $C(2m-5).$ There also exists a semigroup presentation that is strongly $C(3)$ whose induced enveloping group presentation is not $C(2).$
    \label{C(m) implies C(2m-6) sharp}
\end{prop}

\begin{proof} $C(3)$ Case: Suppose $S$ has relations of the following form
\begin{eqnarray*}
    \alpha a\beta &=& \alpha b\beta, \\
    \gamma a \delta &=& \gamma cb\delta,
\end{eqnarray*}
where each of these are non-empty words selected such that $S$ satisfies $C(3)$. (For example, $\alpha, \beta, \gamma, \delta, a, b, c$ could be distinct generators). Then, the induced enveloping group presentation contains the following relators:
\begin{eqnarray*}
    r_1 &=& ab^{-1} \\
   r_2 &=& ab^{-1}c^{-1}.
\end{eqnarray*}
Then, $r_1$ is a proper subword of $r_2$, so the induced enveloping group presentation is not $C(2).$

\medskip

$C(4)$ Case: Suppose $S$ has relations of the following form \begin{eqnarray*}
    \alpha bc\beta &=& \alpha ad\beta \\
    \gamma be\delta&  = &\gamma af \delta \\
    \pi ec\rho &= &\pi fd\rho, 
\end{eqnarray*}
where each of these are non-empty words selected such that $S$ satisfies $C(4)$. (Again one could take $\alpha, \beta, \gamma, \delta, \pi, \rho, a, b, c,d, e, f$ to be distinct generators). The induced enveloping group presentation thus contains the following relators: \begin{eqnarray*}
   r_1&=& bcd^{-1}a^{-1}  \\
    r_2 &=&bef^{-1}a^{-1}\\
    r_3 &=& ecd^{-1}f^{-1}.
\end{eqnarray*}
The cyclic permutation $cd^{-1}a^{-1}b$ of $r_1$ can be decomposed into two group pieces: $cd^{-1}$ which the prefix of the cylic permutation $cd^{-1}f^{1}e$ of $r_3$; and $ab^{-1}$ which is the prefix of the cyclic permutation $a^{-1}bef^{-1}$ of $r_2$. Thus, the induced enveloping group presentation presentation is not $C(3).$

\medskip

$C(5)$ Case: Suppose $S$ has relations of the form 
\begin{eqnarray*}
    \alpha bpc\beta &=& \alpha apd\beta \\
    \gamma b   qc\delta&  = &\gamma aqd\delta,
\end{eqnarray*}
where each of these are non-empty words selected such that $S$ satisfies $C(5)$. Then, the induced enveloping group presentation has the following relators \begin{eqnarray*}
   r_1&=& bpcd^{-1}p^{-1}a^{-1}  \\
    r_2 &=&bqcd^{-1}q^{-1}a^{-1}.
\end{eqnarray*}
The cyclic permutation $p(cd^{-1})p^{-1}(a^{-1}b)$ of $r_1$ can be decomposed into four group pieces, so the induced enveloping group presentation is not $C(5).$

\medskip
$C(m)$ Case for $m\geq 6$: Let $m\geq 6$ and suppose $S$ has relations of the form
\begin{eqnarray*}
    \alpha bp_1\cdots p_{m-4}c\beta &=& \alpha aq_1\cdots q_{m-4} d\beta \\
    \alpha bp_{m-4}\cdots p_1c\beta &=& \alpha aq_{m-4}\cdots q_1 d\beta 
\end{eqnarray*}
so that each relation word is the product of at most $m$ semigroup pieces. Then, the induced enveloping group presentation has the following relators
\begin{eqnarray*}
    r_1 &=& bp_1\cdots p_{m-4}cd^{-1}q^{-1}_{m-4}\cdots q^{-1}_1a \\
    r_2 &=& bp_{m-4}\cdots p_1cd^{-1}q^{-1}_{1}\cdots q^{-1}_{m-4}a
\end{eqnarray*}
The cyclic permutation $$p_1\cdots p_{m-4}(cd^{-1})q^{-1}_{m-4}\cdots q_1^{-1}(a^{-1}b)$$ of $r_1$ can be decomposed into $2m-6$ group pieces, so the induced enveloping group presentation is not $C(2m-5).$

\end{proof}

Thus the following sets of conditions on a semigroup presentation ensure that the induced group presentation satisfies the $C(7)$ condition, sufficient to imply group hyperbolicity: strongly $C(7)$; strongly $C(5)$ and $D(2)$; or strongly $C(4)$ and $D(3).$ Moreover, this result cannot be weakened by appealing to the alternative $C(m)$--$T(q)$ hyperbolicity criteria, since even arbitrarily strong semigroup small overlap conditions do not guarantee that the induced group presentation satisfies $T(4)$, see the following as an example.

\begin{exmpl}
     Consider the following strongly $C(3)$--$D(3)$ semigroup presentation $$S = \langle a, b, c, d \mid abb = bcc, bdd = caa, cba = adb\rangle.$$ The relators of the induced enveloping group presentation are 
\begin{eqnarray*}
    r_1 &= & abbc^{-1}c^{-1}b^{-1} \\
    r_2 &=& bdda^{-1}a^{-1}c^{-1} \\
    r_3 &=& cbab^{-1}d^{-1}a^{-1}
\end{eqnarray*}
These violate $T(4)$, since each of the cyclically adjacent products $r_1r_2$, $r_2r_3$, and $r_3r_1$ admits free cancellation, and no adjacent pair consists of mutually inverse relators.
\end{exmpl}

Thus, even among strongly \(C(m)\) semigroup presentations, the small overlap conditions sufficient to establish non-amenability directly are substantially weaker than those required to guarantee that the induced enveloping group presentation satisfies the classical small cancellation criteria for hyperbolicity. In particular, the result in Corollary~\ref{thm: C(3) non-amenable} cannot be recovered by applying these criteria to the induced enveloping group presentation.

\end{document}